\documentclass[12pt]{amsart}
\usepackage{latexsym, amsmath, amssymb,amsthm,amsopn,amsfonts}
\usepackage{version}
\usepackage{epsfig,graphics,color,graphicx,graphpap}
\usepackage{amssymb}
\usepackage{todonotes}
\usepackage{listings}
\usepackage{mathrsfs}
\usepackage[citecolor=blue]{hyperref}
\usepackage[framemethod=tikz]{mdframed}

\newcommand{\redsout}{\bgroup\markoverwith{\textcolor{red}{\rule[0.5ex]{2pt}{.4pt}}}\ULon}
\numberwithin{equation}{section}

\newcommand{\curl}{\operatorname{curl}}

\newcommand{\p}{\partial}
\newcommand{\C}{\mathbb{C}}
\newcommand{\R}{\mathbb{R}}

\newcommand{\supp}{\text{supp}}

\newcommand{\LV}{\left|}
\newcommand{\RV}{\right|}
\newcommand{\LC}{\left(}
\newcommand{\RC}{\right)}

\definecolor{skyblue}{rgb}{0.85,0.85,1}

\newtheorem{theorem}{Theorem}
\newtheorem{proposition}{Proposition}
\newtheorem{corollary}{Corollary}
\newtheorem{lemma}{Lemma}
\newtheorem{definition}{Definition}

\newtheorem{remark}{Remark}[section]

\title[Non-linear Fractional Magnetic Schr\"odinger Equation]{An Inverse Problem for Non-linear Fractional Magnetic Schr\"odinger Equation} 
 
\begin{document}

\author[Lai]{Ru-Yu Lai}
\address{School of Mathematics, University of Minnesota, Minneapolis, MN 55455, USA}
\curraddr{}
\email{rylai@umn.edu }

\author[Zhou]{Ting Zhou}
\address{Department of Mathematics, Northeastern University, Boston, MA 02115, USA}
\curraddr{}
\email{t.zhou@northeastern.edu}
\maketitle

\begin{abstract}
In this paper, we study forward problem and inverse problem for the fractional magnetic Schr\"odinger equation with nonlinear electric potential. 
We first investigate the maximum principle for the linearized equation and apply it to show that the problem is well-posed under suitable assumptions on the exterior data. Moreover, we explore uniqueness of recovery of both magnetic and electric potentials.
\end{abstract}

\section{Introduction}
 
We study inverse problem for a nonlinear fractional magnetic Schr\"odinger equation (FMSE)
\begin{equation}\label{eqn:LFMSE}
\left\{\begin{split}
 (-\Delta)^s_Au + a(x,u)=0&\qquad\textrm{ in }\Omega\\
 u=g&\qquad\textrm{ in }\Omega_e:=\R^n\backslash \overline\Omega,
\end{split}\right.\end{equation}
where $\Omega\subset \R^n, n\geq2$ is a bounded domain with smooth boundary, $s\in (0,1)$, 
$A$
represents the magnetic vector potential and $a$
is the nonlinear electric potential. 
Here the operator $(-\Delta)^s_A$ is defined by $(\nabla^s+A)^2$ in Section~\ref{sec:notation} following the model introduced in \cite{Covi2020}.

The nonlocal magnetic Schr\"odinger operator $(-\Delta)_A^s$ is shown to arise as continuous limits of long jump random walks with weights (see \cite{Covi2020}), which extends the classical diffusion process modeled by the Laplacian $-\Delta$ and the nonlocal diffusion phenomena modeled by the fractional Laplacian $(-\Delta)^s$. The fractional Laplacian $(-\Delta)^s$ is also seen arising in stochastic theory as the operators associated with symmetric $\alpha$-stable L\'evy processes, for example, in a pricing model in financial mathematics (see \cite{ContTankovbook}). The equations with subcritical nonlinearities were studied in \cite{BRSbook} from the variational point of view.

The inverse problem considered in this paper is a natural generalization of the problems for linear nonlocal fractional elliptic equations such as the fractional Schr\"odinger operator $(-\Delta)^s+q$ and the fractional magnetic Sch\"odinger operator $(-\Delta)^s_A+q$. The fractional analogue of Calder\'on problem for the former operator was first considered in \cite{ghosh2016calder}. More results for Calder\'on problems associated to linear fractional Laplacian can be found in \cite{BGU18, cekic2020calderon, ghosh2017calder, GRSU18, harrach2017nonlocal-monotonicity, harrach2020monotonicity, LLR2019calder, LaiOhm20, ruland2017fractional} and so on. For the fractional magnetic Schr\"odinger operator $(-\Delta)_A^s+q(x)$, which is a nonlocal modification of the classical magnetic Schr\"odinger opertor $(-i\nabla+A(x))^2+q(x)$, the inverse problems were studied in \cite{Li2019} and \cite{Covi2020}. The definition of $(-\Delta)_A^s$ in \cite{Li2019} is given by 
\[(-\Delta)_A^s u(x):=C_{n,s}\textrm{p.v.}\int_{\R^n}\frac{u(x)-e^{i(x-y)\cdot A(\frac{x+y}2)}u(y)}{|x-y|^{n+2s}}~dy,\]
which was first introduced in \cite{PietroSquassina2018}. Based on this definition, the inverse problem for the semilinear equation was considered in \cite{Li2020}. Our paper focuses on the semilinear equation based on the definition of the magnetic Schr\"odinger operator introduced in \cite{Covi2020}. The inverse problems for the semilinear fractional Schr\"odinger equation, where the magnetic potential $A=0$, were also studied in \cite{lai2019global,LaiL2020}.

Meanwhile, recent development of the (higher order) linearization approach for inverse problems of nonlinear elliptic differential equations in order to recover more medium properties also motivates the study of the similar problems for nonlocal diffusion processes.  The inverse boundary value problem for the linear magnetic Schr\"odinger equation $(-i\nabla+A)^2u+qu=0$, has been considered in \cite{Chung2014, ER95, Dos2009, GT2011, Haberman16, IUY2012, KU14, Lai2011, NSU95, PSU10}. 
Specifically, due to a gauge invariance, one can only expect to recover uniquely the magnetic field $\curl A$ and $q$ from the boundary Dirichlet-to-Neumann (DN) map. In dealing with the inverse problems for nonlinear PDEs, a standard approach based on the first order linearization of the DN-map was introduced to identify the linear reaction from the medium, then the full nonlinear medium for certain cases.
See for instance \cite{Sun2002, Isakov93, victor01, victorN, isakov1994global, SunUhlmann97} for the demonstration of the approach in solving the inverse problems for certain semilinear, quasilinear elliptic equations and parabolic equations. 
Recently the higher order linearization of the DN-map has been applied in determining the full nonlinearity of the medium for several different equations. 
The method was successfully applied to solve several inverse problems for nonlinear hyperbolic equations on the spacetime \cite{KLU2018}, where in contrast the underlying problems for linear hyperbolic equations are still open, see also \cite{CLOP,LUW2018} and the references therein. In particular, the second order linearization of the nonlinear boundary map was studied in \cite{CNV2019, Kang2002, Sun96, SunUhlmann97} for nonlinear elliptic equations. Moreover, this higher order linearization technique was also applied to study elliptic equations with power-type nonlinearities, see \cite{FO2019, KU201909, KU201905, LaiL2020, LLLS201905, LLLS201903, Lin202004}. A demonstration of the method can be found in \cite{AZ2018, AZ2020} on nonlinear Maxwell's equations, in \cite{LUY2020} on nonlinear kinetic equations, and in \cite{LLPT2020_1} on semilinear wave equations.   In \cite{LaiZhou2020}, we solved an inverse problem for the magnetic Schr\"odinger equation with nonlinearity in both magnetic and electric potentials $A$ and $q$. \\

In order to describe our main results, we first introduce some definitions, notations and assumptions on the potentials. 
\begin{definition}\label{def:A}
Let $A\in C^s_c (\R^n\times\R^n,\R^n)$. We define the symmetric $A_s$, antisymmetric $A_a$, parallel $A_{||}$ and perpendicular $A_\perp$ parts of $A$ at point $(x,y)$ as
\begin{eqnarray*}&A_s(x,y):=\frac{A(x,y)+A(y,x)}{2},\qquad A_a(x,y):=A(x,y)-A_s(x,y)=\frac{A(x,y)-A(y,x)}2,\\
& A_{||}(x,y):=\left\{\begin{array}{ll}\frac{A(x,y)\cdot(x-y)}{|x-y|^2}(x-y)\quad&\textrm{ if } x\neq y\\
A(x,y)\quad&\textrm{ if }x=y
\end{array}\right. ,\qquad A_{\perp}(x,y):=A(x,y)-A_{||}(x,y).
\end{eqnarray*}
\end{definition}
For the inverse problem, we assume that $A(x,y)\in C^s_c (\R^n\times\R^n,\R^n)$ has compact support in $\Omega\times\Omega$ and satisfies
\begin{equation}\label{eqn:MP_cond_A1}
{ A_{s||}\in H^s(\R^{2n}),\qquad (\nabla\cdot)^sA_{s||}\in L^\infty(\R^n),\qquad }A_{a||}(x,y)\cdot(y-x)\geq0\qquad \textrm{ in }\;\R^n\times\R^n;
\end{equation}
$a(x,z):~\overline\Omega\times\R\rightarrow{\R}$ satisfies
\begin{align}\label{condition a}
\begin{cases}
a(x,0)=0\qquad \hbox{for all }x\in\overline\Omega, \\
\hbox{the map $z\mapsto a(\cdot,z)$ is analytic with values in {$C^s(\overline\Omega)$}}.
\end{cases}
\end{align}
{Therefore, the potential $a$ admits the following Taylor expansion 
\[a(x,z)=\sum_{k=1}^\infty \p_z^ka(x,0)\frac{z^k}{k!},\qquad \p_z^ka(x,0)\in C^s(\overline\Omega).\]
and the convergence of this series is in $C^s(\overline\Omega)$ topology.
}
Together they also satisfy
{
\begin{equation}\label{eqn:MP_cond_A2}
(\nabla\cdot)^sA_{s||}+\int_{\R^n}|A(x,y)|^2~dy+\p_za(x,z)\geq0\qquad \textrm{ for } x\in {\Omega},\; |z|<R_0
\end{equation}
for some constant $R_0>0$.
}
Here we denote by $C^s(\overline\Omega)$ the usual H\"older space.

In Theorem~\ref{THM:well-posedness}, it is shown that there exists a small constant $\varepsilon_0>0$ and a constant $C>0$ such that when the exterior data 
$g \in \mathcal{S}_{\varepsilon_0}$, denoted by 
\[\mathcal{S}_{\varepsilon_0}:=\{g\in C^\infty_c(\Omega_e):\, \|g\|_{C^\infty_c(\Omega_e)}\leq \varepsilon_0\},\]
the problem \eqref{eqn:LFMSE} has a small unique solution $u\in {C^s(\R^n)}$ satisfying $\|u\|_{C^s(\R^n)}\leq C\|g\|_{C^\infty_c(\Omega_e)}$. Therefore, we can  define the Dirichlet-to-Neumann (DN) map
{$
\Lambda^s_{A,a}
$
via the bilinear form as in Lemma \ref{Lemma:DN-map}.}
\\
{We introduce the {\it gauge equivalence} $\sim$ defined in \cite{Covi2020}. We say that two pairs of coefficients $(A_1(x,y),q_1(x))$ and $(A_2(x,y),q_2(x))$ satisfy $(A_1, q_1)\sim(A_2, q_2)$ if and only if  
	{
$$
(-\Delta)^s_{A_1}u + q_1u=(-\Delta)^s_{A_2}u + q_2u
$$	
for all $u\in H^s(\R^n)$. \cite[Lemma 3.8]{Covi2020} further implies that $(A_1, q_1)\sim(A_2, q_2)$ holds if and only if 
	}
\[A_{1,a||}(x,y)=A_{2,a||}(x,y)\quad \textrm{ in } \R^n\times \R^n ,\]
and
\[\quad \int_{\R^n}|A_1|^2~dy+(\nabla\cdot)^sA_{1,s||}+q_1 =\int_{\R^n}|A_2|^2~dy+(\nabla\cdot)^sA_{2,s||}+q_2 \quad \textrm{ in } \Omega .\]

Our main result is stated here. }
\begin{theorem}\label{THH:uniqueness} 
	 	Let $0<s<1$ and let $\Omega \subset \R^n$, $n\geq {2}$ be a bounded domain with smooth boundary. Let $W_1$ and $W_2$ be two arbitrary nonempty open subsets in $\Omega_e$. Suppose that $(A_j, a_j)$ satisfy \eqref{eqn:MP_cond_A1}, \eqref{condition a} and \eqref{eqn:MP_cond_A2} for $j=1,2$. 
Then if 
	 \begin{align}\label{DN map in Thm 1}
	 \left.	\Lambda^s_{A_1,a_1}[g] \right|_{W_2} = 	  \left.	\Lambda^s_{A_2,a_2}[g] \right|_{W_2} \qquad \text{ for any }g\in {\mathcal S_{\varepsilon_0}\cap C_c^\infty(W_1)},
	 \end{align} 
where $\varepsilon_0>0$ is sufficiently small $($Here $[g]$ is the equivalence class of $g$ in $H^s(\R^n)\backslash\widetilde H^s(\Omega)$$)$, we have 
	 \begin{align}\label{first order}
	 (A_1,\p_za_1(x,0))\sim(A_2,\p_za_2(x,0)) 
	 \end{align}
	 and 
	 \begin{align}\label{higher order}
		 a_1(x,z)-\p_za_1(x,0)z = a_2(x,z)-\p_za_2(x,0)z\quad \textrm{ in }\Omega\times \mathbb{R}.
	 \end{align}
\end{theorem}
\vskip.2cm

	The proof of Theorem~\ref{THH:uniqueness} is built upon several preliminary results: the well-posedness for the problem, the Runge approximation property and maximum principle. We start by investigating the forward problem since the study of the inverse problem stands on it. 
	To this end, we formulate the maximum principle and a barrier function, which can be applied to prove the boundedness of solution to the FMSE. Together with the benefit introduced by the nonlinearity of the equation, they guarantee the effectiveness of the fixed point theorem. This then leads to the well-posedness result for the nonlinear equation under study. 
	Moreover, to reconstruct unknown potentials, we apply the higher order linearization scheme. By taking derivatives of the integral identity for the DN map multiple times, the Runge approximation property, which states that the set of solutions is dense in $L^2(\Omega)$, then plays in to extract the information of unknown potentials out of the integral.

This theorem guarantees uniqueness of higher order term of $a(x,z)$. However, similar to the inverse problem for the linear FMSE, it is expected that one can only determine the magnetic and linear electric potentials up to a gauge as shown in \eqref{first order}. In fact, in Section~\ref{sec:inverse}, we found that only coefficients $A$ and $\p_za(x,0)$ appear in the linearized equation and then they are recovered up to the gauge by applying the available inverse problem result for the linear equation. 
While in the higher order linearization steps, the higher order coefficients $\p_z^ka(x,0),k\geq 2$ can be viewed as a part of the source term in the $k$-th linearized equation. Therefore, this explains the unique reconstruction of these terms without any gauge in \eqref{higher order}.

\vskip.2cm
The structure of the paper is as follows. In Section~\ref{sec:forward}, we establish the well-posedness of the problem \eqref{eqn:LFMSE} for a small enough exterior data $g$ by deriving the maximum principle and the barrier function in order to obtain the $C^s$ regularity of the solution for $(-\Delta)^s_A+q$. 
In Section~\ref{sec:inverse}, we determine the potentials $A$ and $a(x,z)$ using the linearization steps.

\section{The forward problem}\label{sec:forward}

\subsection{Notations}\label{sec:notation}
We introduce the notations and properties below.  
We define $H^s(\R^n)=W^{s,2}(\R^n)$ to be the standard $L^2$-based Sobolev space with norm 
$$
    \|u\|_{H^s(\R^n)} = \|\mathcal{F}^{-1}((1+|\xi|^2)^{s/2}\mathcal{F} u)\|_{L^2(\R^n)},
$$
where $\mathcal{F}$ is the Fourier transform defined as 
$$
    \mathcal{F} u (\xi)= \int_{\R^n} e^{-ix\cdot\xi} u(x)\,dx.
$$

Let $U$ be an open set in $\R^n$. For scalar $\beta\in\mathbb{R}$, we define the following spaces:
\begin{align*}
H^{\beta}(U) & :=\left\{u|_{U}:\, u\in H^{\beta}(\mathbb{R}^{n})\right\},\\
\widetilde{H}^{\beta}(U) & :=\text{closure of \ensuremath{C_{c}^{\infty}(U)} in \ensuremath{H^{\beta}(\mathbb{R}^{n})}},\\
H_{0}^{\beta}(U) & :=\text{closure of \ensuremath{C_{c}^{\infty}(U)} in \ensuremath{H^{\beta}(U)}}.
\end{align*}

Following the notations in \cite{Covi2020}, the magnetic fractional Laplacian $(-\Delta)^s_A$ is an operator mapping from $H^s(\R^n)$ to $H^{-s}(\R^n)$, such that for all $u, v\in H^s(\R^n)$, 
\[\langle (-\Delta)^s_A u, v\rangle =\langle \nabla^s_A u, \nabla^s_A v\rangle. \]  
The magnetic fractional gradient operator $\nabla_A^s$ is defined by $\nabla^s+A(x,y)$ with the fractional gradient $\nabla^s: H^s(\R^n)\rightarrow L^2(\R^{2n})$ extends the definition
\[\nabla^s u(x,y)=\sqrt{\frac{C_{n,s}}{2}}\frac{u(x)-u(y)}{|y-x|^{n/2+s+1}}(y-x).\]
Then the fractional divergence $(\nabla\cdot)^s: L^2(\R^{2n})\rightarrow H^{-s}(\R^n)$ is defined by 
\[\langle (\nabla\cdot)^su,v\rangle=\langle u,\nabla^s v\rangle\qquad\textrm{ for } v\in H^s(\R^n).\]
Recall that the fractional Laplacian $(-\Delta)^s: H^s(\R^n)\rightarrow H^{-s}(\R^n)$ is defined by
\[(-\Delta)^s u(x):=C_{n,s}\textrm{ p.v.} \int_{\R^n}\frac{u(x)-u(y)}{|x-y|^{n+2s}}~dy ,\]
where the constant $C_{n,s}$ depends only on $n,s$, see for instance \cite{di2012hitchhiks}, and p.v. stands for the principal value.
Then we have $(-\Delta)^s=(\nabla\cdot)^s\nabla^s$ in weak sense, that is, $\langle (-\Delta)^s u, v\rangle=\langle \nabla^s u,\nabla^s v\rangle$ 
for $u, v\in H^s(\R^n)$.

\subsection{Preliminary results}\label{sec:preliminary} 
The proof of the following result can be found in \cite[Lemma~3.15]{Covi2020}, where the regularity of $A$ can be relaxed if certain integrability conditions are imposed.  
\begin{proposition}[The Runge approximation] \label{prop:RAP}Suppose that $q\in L^p(\Omega)$, $p:=\max\{2,n/(2s)\}$ and $A \in  C^s_c(\R^n\times\R^n,\C^n)$ has compact support in $\Omega\times\Omega$ and satisfies \eqref{eqn:MP_cond_A3}.
	Let $W$ be an open set in $\Omega_e$ and $u_g$ be the solution to $(-\Delta)^s_Au_g+qu_g=0$ in $\Omega$ with $u_g	=g$ in $\Omega_e$.
	Then the set $\mathcal{R}=\{u_g|_\Omega:\, g\in C_c^\infty(W)\}$ is dense in $L^2(\Omega)$.
\end{proposition}

The Runge approximation will be applied to recover the unknowns $A$ and $a(x,z)$ in Section~\ref{sec:inverse}.

\subsection{Boundedness of solutions}
We will follow the steps in \cite{lai2019global} to prove an $L^\infty$-bound for the weak solution of \eqref{eqn:FMSE_linear_fg} in Proposition~\ref{prop:L-infty}. To this end, we prove a maximum principle in Proposition~\ref{prop:MP_LFMSE} and construct a barrier function in Lemma~\ref{lem:barrier}.

\begin{proposition}[Maximum principle for the weak solution of FMSE]\label{prop:MP_LFMSE}
Let $\Omega$ be a bounded Lipschitz domain in $\R^n$ and $A(x,y)\in C^s_c(\R^n\times\R^n,\R^n)$ satisfy  
\begin{equation}\label{eqn:MP_cond_A3}
{A_{s||}\in H^s(\R^{2n}) ,\qquad }A_{a||}(x,y)\cdot(y-x)\geq0\qquad \textrm{ in }\;\R^n\times\R^n;
\end{equation}
and 
\begin{equation}\label{eqn:MP_cond_A4}
(\nabla\cdot)^sA_{s||}+\int_{\R^n}|A(x,y)|^2~dy\geq0\qquad \textrm{ for } x\in {\Omega}.
\end{equation}
Suppose $u\in H^s(\Omega)$ is a weak solution of
\[\left\{\begin{split}
(-\Delta)^s_Au=F\qquad&\textrm{ in }\;\Omega,\\
u=g\qquad & \textrm{ in }\;\Omega_e.
\end{split}\right.\]
Then if $0\leq F\in L^\infty(\Omega)$ and $0\leq g\in L^\infty(\Omega_e)$, we have $u\geq0$ in $\Omega$, hence in $\R^n$. 
\end{proposition}

\begin{proof}
Using the weak formulation, we obtain for $\phi\in H^s_0(\Omega)$ satisfying $\phi\geq0$, we have
\begin{align}\label{prop:phi}
\langle \nabla^su(x)+A(x,y)u(x), \nabla^s\phi(x)+A(x,y)\phi(x)\rangle=\int_\Omega F\phi~dx\geq0,
\end{align}
where we used $F\geq 0$ in $\Omega$.
Let $u^-:=\max\{-u,0\}$. Then by $u\in H^s(\Omega)$, we can take
$\phi=u^-\in H^s_0(\Omega)$ as a test function. 
We want to show that if $\phi\not\equiv 0$, the left hand side of \eqref{prop:phi} turns out to be negative, in order to draw a contradiction. We rewrite the left hand side of \eqref{prop:phi} as
\[\begin{split}
\langle \nabla^su+Au,\nabla^s\phi+A\phi\rangle&=\langle\nabla^su,\nabla^s\phi\rangle+\left[\langle\nabla^su, A(x,y)\phi(x)\rangle-\langle\nabla^su, A(y,x)\phi(x)\rangle\right]\\
&\quad +\left[\langle\nabla^su,A(y,x)\phi(x)\rangle+\langle\nabla^s\phi,A(x,y)u(x)\rangle+\langle Au, A\phi\rangle\right]\\
&=: \langle\nabla^su,\nabla^s\phi\rangle+I+II.
\end{split}\]
It was shown in \cite{lai2019global} (see the proof of Proposition 3.1 in \cite{lai2019global}) that the fraction Laplacian term $\langle\nabla^su,\nabla^s\phi\rangle<0$, where $g\geq0$ in $\Omega_e$ is used. We will discuss below that $I+II$ is indeed nonpositive, that is $I+II\leq 0$, which then leads to a contradiction to \eqref{prop:phi}.

To this end, the term $I$ is actually $2\langle\nabla^su, 	A_a(x,y)\phi(x)\rangle$ and 
\[
I=2\sqrt{\frac{C_{n,s}}2}\int_{\R^{2n}}\frac{A_a(x,y)\cdot(y-x)}{|x-y|^{n/2+s+1}}\left(u(x)-u(y)\right)\phi(x)~dxdy.
\]
Note that the integrand function vanishes on $S_0:=\{(x,y)~|~u(x)\geq0\}$ since $\phi=0$ on $S_0$. On the set $S_1:=\{(x,y)~|~u(x)<0\;\textrm{ and }\;u(y)\geq0\},$ we have 
$$(u(x)-u(y))\phi(x)<0\qquad \hbox{ in }S_1,$$ 
by \eqref{eqn:MP_cond_A3} (which is equivalent to $A_a(x,y)\cdot(y-x)\geq0$), and thus the integrand function is $\leq0$.  
On the remaining set $S_2:=\{(x,y)~|~u(x)<0\;\textrm{ and }\;u(y)<0\}$, we have $$\phi(x)=u^-=-u(x)\qquad \hbox{ in }S_2.$$ 
We denote $\Omega^-:=\{x~|~u(x)<0\}$. Then it is sufficient to consider the integral over $\Omega^-\times\Omega^-$:
\[\begin{split}
&-2\sqrt{\frac{C_{n,s}}2}\int_{\Omega^-\times\Omega^-}\frac{A_a(x,y)\cdot(y-x)}{|x-y|^{n/2+s+1}}\left(u(x)-u(y)\right)u(x)~dxdy\\
=&-\sqrt{\frac{C_{n,s}}2}\int_{\Omega^-\times\Omega^-}\frac{A_a(x,y)\cdot(y-x)}{|x-y|^{n/2+s+1}}\left(u(x)-u(y)\right)u(x)\\
&\qquad\qquad\qquad+\frac{A_a(y,x)\cdot(x-y)}{|x-y|^{n/2+s+1}}\left(u(y)-u(x)\right)u(y)~dxdy\\
=&-\sqrt{\frac{C_{n,s}}2}\int_{\Omega^-\times\Omega^-}\frac{A_a(x,y)\cdot(y-x)}{|x-y|^{n/2+s+1}}\left(u(x)-u(y)\right)^2~dxdy \leq 0,
\end{split}\]
since $A_a(y,x)=-A_a(x,y)$ and also by \eqref{eqn:MP_cond_A3}. Combining these estimates together, we have $I\leq 0$. 

Lastly, following the proof of Lemma 3.3 in \cite{Covi2020}, the term $II$ satisfies 
\[II= \int_{\Omega}\left((\nabla\cdot)^sA_{s||}(x)+\int_{\R^n}|A(x,y)|^2~dy\right)\phi(x)u(x)~dx\leq 0,\]
by the assumption \eqref{eqn:MP_cond_A4} and the fact $u(x)\phi(x)\leq0$ in $\R^n$. 

This completes the proof of $\langle \nabla^su+Au,\nabla^s\phi+A\phi\rangle<0$, which in turn concludes the proof of the proposition by contradiction.
\end{proof}

We now build a barrier function.

\begin{lemma}[Barrier]\label{lem:barrier}
Let $\Omega$ be a bounded Lipschitz domain in $\R^n$. Assume that $A\in C^s_c (\R^n\times\R^n;\R^n)$ satisfies {\eqref{eqn:MP_cond_A3} and }\eqref{eqn:MP_cond_A4}. Then there exists $\varphi\in C_c^\infty(\R^n)$ satisfying 
\[\left\{\begin{split}
(-\Delta)^s_A \varphi \geq1&\qquad \textrm{ in }\;\Omega,\\
\varphi\geq 0 &\qquad \textrm{ in }\;\R^n,\\
\varphi\leq C &\qquad \textrm{ in }\;\Omega,
\end{split}\right.\]
for some constant $C$ depending on $n, s$ and $\Omega$.
\end{lemma}
\begin{proof}
We will show the barrier function $\varphi$ in Lemma 3.4 in \cite{lai2019global} qualifies. More specifically, let $B_R$ be a large ball such that $\Omega\Subset B_R$ and $\eta\in C^\infty_c(B_R)$ be a smooth cutoff function satisfying
\[0\leq\eta\leq1\quad\textrm{ in }\;\R^n,\qquad \eta\equiv 1\quad\textrm{ in }\;\Omega.\]
We directly use the format of $(-\Delta)^s_A$ given in Lemma 3.3 of \cite{Covi2020},
\begin{align}\label{Id fractional M}
(-\Delta)^s_A\eta=(-\Delta)^s\eta+2\int_{\R^n}\left(A_{a||}\cdot\nabla^s\eta\right)~dy+\left((\nabla\cdot)^sA_{s||}+\int_{\R^n}|A|^2~dy\right)\eta.
\end{align}
By the definition of $(-\Delta)^s$ and the fact that $\eta$ has maximum value $1$ in $\Omega$, one has when $x\in\Omega$
\[\begin{split}(-\Delta)^s\eta (x)
=&\,C_{n,s}\int_{\R^n}\frac{\eta(x)-\eta(z)}{|x-z|^{n+2s}}~dz\\
\geq&\, C_{n,s} \int_{\R^n\backslash B_R}\frac{1}{|x-z|^{n+2s}}~dz\\ 
\geq&\,\frac{C_{n,s}}2\int_{\R^n\backslash B_R}\frac{1}{(R+|z|)^{n+2s}}~dz:=\lambda>0,
\end{split}\]
where $\lambda$ depends only on $n, s$ and $\Omega$. By the condition \eqref{eqn:MP_cond_A4} and $A$ is compactly supported in $\Omega\times\Omega$, we only need to show 
$$2\int_{\R^n}\left(A_{a||}\cdot\nabla^s\eta\right)~dy=2\int_{\Omega}\left(A_{a||}\cdot\nabla^s\eta\right)~dy\geq 0$$
when $x\in\Omega$.
This is verified because the integrand
\[
\left(A_{a||}\cdot\nabla^s\eta\right)(x,y)=\sqrt{C_{n,s}\over 2}\frac{A_a(x,y)\cdot(y-x)}{|x-y|^{n/2+s+1}}(\eta(x)-\eta(y))\geq 0 \qquad \hbox{ when $x\in\Omega$.}
\]
Here we applied \eqref{eqn:MP_cond_A3} to obtain $A_a(x,y)\cdot(y-x)=A_{a||}(x,y)\cdot(y-x)\geq 0$ and also observed that $\eta(x)-\eta(y)=1-\eta(y)\geq 0$ for $x\in\Omega$ and $y\in\R^n$. Thus we have shown that $(-\Delta)^s_A\eta \geq \lambda$.
Finally, set $\varphi(x)=\frac{\eta(x)}{\lambda}$. The upper bound $C$ is $1/\lambda$, hence depends on $n, s$ and $\Omega$.
\end{proof}

With the maximum principle and the barrier function, we can show an $L^\infty$-estimate for the weak solution.
\begin{proposition}[$L^\infty$-bound for the weak solution of FMSE]\label{prop:L-infty}
Let $\Omega$ be a bounded smooth domain and $A\in C^s_c (\R^n\times\R^n,\R^n)$ satisfy \eqref{eqn:MP_cond_A3} and \eqref{eqn:MP_cond_A4}. For $F\in L^\infty(\Omega)$ and $g\in L^\infty(\Omega_e)$, assume that $u\in H^s(\Omega)$ is a weak solution of 
\begin{equation}\label{eqn:FMSE_linear_fg}
\left\{\begin{split}
(-\Delta)^s_A u=F&\quad\textrm{ in }\;\Omega,\\
u=g&\quad\textrm{ in }\;\Omega_e.
\end{split}\right.\end{equation}
Then 
\[\|u\|_{L^\infty({\R^n})}\leq \|g\|_{L^\infty(\Omega_e)}+C\|F\|_{L^\infty(\Omega)}\]
for some constant $C>0$ depending on $n, s$ and $\Omega$. (It can be the same constant as in Lemma \ref{lem:barrier}.)
\end{proposition}

\begin{proof}
It is a standard proof such as in \cite{lai2019global}. For completeness we prove it here for $(-\Delta)^s_A$. We set 
\[v(x)=\|g\|_{L^\infty(\Omega_e)}+\|F\|_{L^\infty(\Omega)}\varphi(x)\qquad \hbox{ in }\R^n,\]
where $\varphi$ is the barrier function in Lemma \ref{lem:barrier}. Since $(-\Delta)^s_A \varphi\geq 1$ in $\Omega$, $$(-\Delta)^s_A (\|g\|_{L^\infty(\Omega_e)})= \left((\nabla\cdot)^sA_{s||}+\int_{\R^n}|A|^2~dy\right)\|g\|_{L^\infty(\Omega_e)}$$ and \eqref{eqn:MP_cond_A4}, we have
\[(-\Delta)_A^s v\geq \|F\|_{L^\infty(\Omega)}\geq F=(-\Delta)^s_A u,\]
in $\Omega$. Moreover, we get $v\geq \|g\|_{L^\infty(\Omega_e)}\geq u$ in $\Omega_e$ due to $\varphi\geq 0$ in $\Omega_e$. 
Applying Proposition \ref{prop:MP_LFMSE} we obtain $v-u\geq 0$ in $\R^n$. This proves
\[u(x)\leq \|g\|_{L^\infty(\Omega_e)}+C\|F\|_{L^\infty(\Omega)}\qquad \hbox{ in }\Omega.\]
Similarly, the same would hold for $-u$. This completes the proof. 
\end{proof}

\begin{remark}\label{remark:q}
Here we can certainly add a scalar potential $q\in L^\infty(\Omega,\R)$ to obtain the above results {(Proposition~\ref{prop:MP_LFMSE}, Lemma~\ref{lem:barrier}, Proposition~\ref{prop:L-infty})} for equation $(-\Delta)^s_Au+q(x)u=F$ with the condition \eqref{eqn:MP_cond_A4} replaced by 
\begin{equation}\label{eqn:MP_cond_A5}
(\nabla\cdot)^sA_{s||}+\int_{\R^n}|A(x,y)|^2~dy + q(x)\geq0\qquad \textrm{ for } x\in \Omega.\end{equation}
\end{remark}

We also need the following strong maximum principle in the proof of Theorem \ref{THH:uniqueness}.
\begin{proposition}[Strong maximum principle]\label{prop:MP_Strong}
	Let $\Omega$ be a bounded Lipschitz domain in $\R^n$. {Let $A(x,y)\in C^s_c (\R^n\times\R^n,\R^n)$ and $q\in L^\infty(\Omega,\R)$ satisfy conditions \eqref{eqn:MP_cond_A3} and \eqref{eqn:MP_cond_A5}.} 
Suppose $u\in H^s(\Omega)$ is a weak solution of
\[\left\{\begin{split}
(-\Delta)^s_Au + qu=F\qquad&\textrm{ in }\;\Omega,\\
u=g\qquad & \textrm{ in }\;\Omega_e.
\end{split}\right.\]
Then if $0\leq F\in L^\infty(\Omega)$ and $0\leq g\in L^\infty(\Omega_e)$ with $g \not\equiv 0$, we have $u> 0$ in $\Omega$.
\end{proposition}	
\begin{proof}
	Proposition~\ref{prop:MP_LFMSE} yields that $u\geq 0$ in $\R^n$. Suppose that $u$ is not strictly positive in $\Omega$. {Then there must exists a nonempty subset $B\subset \Omega$ with positive measure so that $u=0$ in $B$. We take a smooth function $0\neq\varphi\in C^\infty_c(B)$ satisfying $\varphi\geq 0$ in $B$. Since $u$ is the weak solution, by \eqref{Id fractional M}, we have
   \begin{align}\label{proof of SMP}
	\begin{split}
	0&\leq  \int_B F(x ) \varphi(x)\, dx=\int_B ((-\Delta)^s_A u(x )+q(x )u(x )) \varphi(x)\,dx \\
	&=\int_B [-  C_{n,s} \int_{\R^n} \frac{ u(y)}{|x-y|^{n+2s}}\, dy - \sqrt{2C_{n,s}} \int_{\R^n} \frac{ A_a(x,y)\cdot(y-x) u(y)}{|x-y|^{n/2+s+1}}\, dy ] \varphi(x)\,dx \\
	& \leq  0,
	\end{split}
	\end{align}
	where the last inequality follows by the fact that $u,\varphi\geq 0$, the constant $C_{n,s}>0$ and also $A_a(x,y)\cdot(y-x)\geq 0$ based on the assumption. 	
	This further yields that the first term in \eqref{proof of SMP} satisfies $$\int_B -  C_{n,s} \int_{\R^n} \frac{ u(y)}{|x-y|^{n+2s}}\, dy \varphi(x)dx=0.$$ 
Now due to $C_{n,s}>0$, $u,\varphi\geq 0$, from the above identity, we can then derive $\int_{\R^n} \frac{ u(y)}{|x-y|^{n+2s}}\, dy=0$, which leads to $u\equiv 0$ in $\R^n$, which contradicts to the assumption that $g\not \equiv0$ in $\Omega_e$.}
\end{proof}

\subsection{Well-posedness of a nonlinear fractional MSE}
With the $L^\infty$-bound, we will show that the solution of the fractional magnetic equation \eqref{eqn:FMSE_linear_fg} indeed has $C^s$ regularity. This regularity is essential to prove the well-posedness result later.

\begin{lemma}[$C^s$-estimate]\label{lem:Cs_est}
Let $\Omega$ be a bounded smooth domain in $\R^n$ ($n\geq2$). Suppose $A\in C_c^s(\R^n\times\R^n,\R^n)$ and $q\in L^\infty(\Omega,\R)$ satisfy conditions \eqref{eqn:MP_cond_A3} and \eqref{eqn:MP_cond_A5}. {Suppose also that $(\nabla\cdot)^sA_{s||}\in L^\infty(\R^n)$}. Given $F\in L^\infty(\Omega)$ and {$g\in C^\infty_c(\Omega_e)$}, if $u\in H^s(\R^n)$ is a weak solution of 
\begin{equation} 
\left\{\begin{split}
(-\Delta)^s_A u+qu=F&\quad\textrm{ in }\;\Omega,\\
u=g&\quad\textrm{ in }\;\Omega_e,
\end{split}\right.\end{equation}
 then we have $u\in C^s(\R^n)$ and $u$ satisfies
\begin{equation}\label{eqn:Cs_estimate}
\|u\|_{C^s(\R^n)}\leq C\left(\|F\|_{L^\infty(\Omega)}+\|g\|_{C^s(\Omega_e)}\right),
\end{equation}
where the constant $C$ depends on $A,q,n,s$, and $\Omega$.
\end{lemma}
\begin{proof}
We first extend $g\in C^\infty_c(\Omega_e)$ to the whole $\R^n$ by zero and denote this extension by $g$ as well.
Then $\widetilde u:=u-g$ satisfies
\begin{align}\label{eqn:FMSE_linear_wt_f}
	\left\{\begin{array}{ll}
	(-\Delta)_A^s \widetilde u +q\widetilde u = F-[(-\Delta)_A^s +q]g=:\widetilde F\in L^\infty(\Omega) & \text{ in }\Omega, \\
	\widetilde u=0 & \text{ in }\Omega_e.
	\end{array}  \right.
	\end{align}
Note that since $n\geq 2$, the change of variables yields that $(-\Delta)_A^s g\in L^\infty(\Omega)$.
By Proposition \ref{prop:L-infty} and Remark~\ref{remark:q}, one has
\[\|\widetilde u\|_{L^\infty(\R^n)}\leq C_{n,s,\Omega}\|\widetilde F\|_{L^\infty(\Omega)},\]
where the constant $C_{n,s,\Omega}$ depends on $n,s,\Omega$. Also, from the format of $(-\Delta)^s_A$ given in Lemma 3.3 of \cite{Covi2020}, we have 
\begin{align}
    \label{eqn:FMSE_linear_wt_f 2}(-\Delta)^s \widetilde u=\widetilde F-H_{\widetilde u},\qquad \widetilde{u}|_{\Omega_e}=0,	
\end{align} 
where 
\[H_{\widetilde u}:=2\int_{\R^n}\left(A_{a||}\cdot\nabla^s \widetilde u\right)~dy+\left((\nabla\cdot)^sA_{s||}+\int_{\R^n}|A|^2~dy+q\right)\widetilde u.\]
Note that since $A$ is $C^s$ and compactly supported in $\Omega\times\Omega$ and $n/2\geq 1>s>0$, we can derive that $\frac{A_a(x,y)\cdot(y-x)}{|x-y|^{n/2+s+1}}$ is integrable in $\R^n$ and thus
\begin{align*}
\LV 2\int_{\R^n}\left(A_{a||}\cdot\nabla^s \widetilde u\right)~dy \RV &= \LV\sqrt{2C_{n,s}}\int_{\R^n}\left(A_{a}(x,y)\cdot (y-x) {\widetilde u(x)-\widetilde u(y) \over |x-y|^{n/2+s+1}}\right)~dy \RV\\
&\leq 2\sqrt{2 C_{n,s}}\|\widetilde u\|_{L^\infty(\Omega)}   \int_{\R^n} { |A_{a}\cdot (y-x)|\over |x-y|^{n/2+s+1}} ~dy  \\
&\leq C_{A,n,s,\Omega}\|\widetilde F\|_{L^\infty(\Omega)};
\end{align*}
{and $\left((\nabla\cdot)^sA_{s||}+\int_{\R^n}|A|^2~dy+q\right)$ is bounded in $\Omega$}. This yields that
\[\left|\left((\nabla\cdot)^sA_{s||}+\int_{\R^n}|A|^2~dy+q\right)\widetilde u\right|\leq C_{A,q,n,s,\Omega}\|\widetilde F\|_{L^\infty(\Omega)}\qquad\textrm{ in }\Omega.\]
Therefore,
\[\|H_{\widetilde u}\|_{L^\infty(\Omega)}\leq C_{A,q,n,s,\Omega}\|\widetilde F\|_{L^\infty(\Omega)}.\]
Applying \cite[Proposition~1.1]{ros2014dirichlet} to \eqref{eqn:FMSE_linear_wt_f 2}, we have $\widetilde u\in C^s(\R^n)$ and
\[\|\widetilde u\|_{C^s(\R^n)}\leq C_{\Omega, s}\|\widetilde F-H_{\widetilde u}\|_{L^\infty(\Omega)}\leq C_{A,q,n,s,\Omega}\|\widetilde F\|_{L^\infty(\Omega)},\]
which proves \eqref{eqn:Cs_estimate}.
\end{proof}

We are ready to prove unique existence of solutions to the nonlinear FMSE. 

\begin{theorem}[Well-posedness for the nonlinear equation]\label{THM:well-posedness}
Let $\Omega$ be a bounded smooth domain and $A\in C^s_c (\R^n\times\R^n,\C^n)$ and $a(x,z)$ satisfy \eqref{eqn:MP_cond_A1}, \eqref{condition a} and \eqref{eqn:MP_cond_A2}. 
Then there exists $\varepsilon>0$ small enough such that when {$g\in C_c^\infty(\Omega_e)$} with $\|g\|_{C^\infty_c(\Omega_e)}<\varepsilon$, the boundary value problem 
\begin{align}\label{eqn:nonlinear_FMSE}
	\begin{cases}
	(-\Delta)^s_A u+a(x,u) =0 & \hbox{ in } \Omega,\\
	u=g  &  \hbox{ in } \Omega_e,\\
	\end{cases} 
	\end{align}
admits a unique solution {$u\in C^s(\R^n){\cap H^s(\R^n)}$} satisfying
\begin{equation}\label{eqn:Cs_est}\|u\|_{C^s(\R^n)}\leq C\|g\|_{C^\infty_c(\Omega_e)},\end{equation}
where $C$ is a constant depending on $A,\p_z a(x,0),n,s$, and $\Omega$.
\end{theorem}
\begin{proof}
Combining Section 3.3 \cite{Covi2020} and Lemma~\ref{lem:Cs_est} yields that there exists a unique solution {$u_0\in C^s(\R^n){\cap H^s(\R^n)}$} to the linear equation
\begin{align}\label{Cs estimate}
	\begin{cases}
	(-\Delta)^s_A u_0+\partial_za(x,0)u_0 =0 & \hbox{ in } \Omega,\\
	u_0=g  &  \hbox{ in } \Omega_e,\\
	\end{cases} 
	\end{align}
such that
\begin{align}\label{estimate u0}
\|u_0\|_{C^s(\R^n)}\leq C\|g\|_{C^\infty_c(\Omega_e)}.
\end{align}

Then looking for a solution of \eqref{eqn:nonlinear_FMSE} is equivalent to solving for $v:=u-u_0$ in
\begin{align}\label{eqn:v}
	\begin{cases}
	(-\Delta)^s_A v+\partial_z a(x,0)v =-\left(a(x,u_0+v)-\partial_za(x,0)(u_0+v)\right) & \hbox{ in } \Omega,\\
	v=0  &  \hbox{ in } \Omega_e.\\
	\end{cases} 
	\end{align}
To achieve this, let us define the set 
	$$X_\delta = \left\{\phi\in C^{s}(\R^n):\ \phi|_{\Omega_e}=0,\ \|\phi\|_{C^{s}(\R^n)}\leq \delta \right\},$$
	where $0 < \delta < 1$ will be determined later. It is easy to see that $X_\delta$ is a Banach space. 	
We also define the map $\mathcal{F}:X_\delta\rightarrow L^\infty(\Omega)$, by
$$\mathcal F(v):=-\left(a(x,u_0+v)-\partial_za(x,0)(u_0+v)\right),$$ and by Section 3.3 \cite{Covi2020} and Lemma~\ref{lem:Cs_est} the operator
$$\mathcal L_s^{-1}: F\in L^\infty(\Omega)~\mapsto~u_{F,0}\in C^s(\R^n){\cap H^s(\R^n)} $$ is bounded, where $u_{F,0}$ denotes the solution of $\left[(-\Delta)^s_A+\partial_za(x,0)\right]u_{F,0}=F$ in $\Omega$ and $u_{F,0}=0$ in $\Omega_e$. It suffices to show that $\mathcal L_s^{-1}\circ\mathcal F$ is a contraction map on $X_\delta$.

To show $\mathcal L_s^{-1}\circ\mathcal F$ is contractive, we first apply Taylor's theorem and $a(x,0)=0$ to obtain
\[a(x,z)=\partial_za(x,0)z+a_r(x,z)z^2,\qquad (x,z)\in \Omega\times \R\]
with
\[a_r(x,z):=\int_0^1\partial_z^2a(x,tz)(1-t)~dt.\]
Therefore,
\[\mathcal F(v)=-a_r(x,u_0+v)(u_0+v)^2.\]
In $\Omega$, it is not hard to see that $a_r(x,u_0+v)$ is bounded by a constant. Then 
\begin{equation}\label{eqn:F_est}\|\mathcal F(v)\|_{L^\infty(\Omega)}\leq C\left(\|u_0\|_{C^s(\R^n)}+\|v\|_{C^s(\R^n)}\right)^2\leq C(\varepsilon+\delta)^2\end{equation}
for $v\in X_\delta$. This implies
\[\|\mathcal L_s^{-1}\circ\mathcal F(v)\|_{C^s(\R^n)}\leq C(\varepsilon+\delta)^2\leq \delta\]
when {$\varepsilon<\mathcal C\delta$ for some $\mathcal C>0$ (this is to say that $\delta$ cannot be arbitrarily small, but depending on $\varepsilon$}) and $\delta$ is small enough, that is, $\mathcal L_s^{-1}\circ\mathcal F$ maps $X_\delta$ to itself.

To show it is a contraction, we derive for $v_1, v_2\in X_\delta$
\[\begin{split}\|\mathcal F(v_1)-\mathcal F(v_2)\|_{L^\infty(\Omega)}&\leq    \|a_r(x,u_0+v_1)-a_r(x,u_0+v_2)\|_{L^\infty(\Omega)}\left(\|u_0\|_{C^s(\R^n)}+\|v_1\|_{C^s(\R^n)}\right)^2\\
&\quad+\|a_r(x,u_0+v_2)\|_{L^\infty(\Omega)}\|2u_0+v_1+v_2\|_{C^s(\R^n)}\|v_1-v_2\|_{C^s(\R^n)}\\
&\leq C[(\varepsilon+\delta)^2+2\varepsilon+2\delta]\|v_1-v_2\|_{C^s(\R^n)},
\end{split}\]
where we used the fact that $a_r(x,z)$ is Lipchitz in $z$. This shows
\[\|\mathcal L_s^{-1}(\mathcal F(v_1)-\mathcal F(v_2))\|_{C^s(\R^n)}\leq C[(\varepsilon+\delta)^2+2\varepsilon+2\delta]\|v_1-v_2\|_{C^s(\R^n)},\]
hence $\mathcal L_s^{-1}\circ\mathcal F$ is a contraction when $\varepsilon$ and $\delta$ are small enough. 

Therefore, the contraction mapping principle yields that there exists a fixed point $v\in X_\delta$ so that $v=\mathcal L_s^{-1}\circ\mathcal F(v)$ satisfies \eqref{eqn:v}.
Lastly, we have by \eqref{estimate u0} and \eqref{eqn:F_est}
\[\|v\|_{C^s(\R^n)}=\|\mathcal L_s^{-1}\circ\mathcal F(v)\|_{C^s(\R^n)}\leq C\|\mathcal F(v)\|_{L^\infty(\Omega)}\leq C(\varepsilon+\delta)(\|g\|_{C^\infty_c(\Omega_e)}+\|v\|_{C^s(\R^n)}),\]
which gives 
$$
\|v\|_{C^s(\R^n)} \leq C \|g\|_{C^\infty_c(\Omega_e)}
$$
if $\varepsilon,\delta$ are sufficiently small.
Combining the above estimate for $v$ with \eqref{estimate u0}, this gives \eqref{eqn:Cs_est}.
\end{proof}

\begin{corollary}\label{Corollary:Hs_reg}
Assume that $\Omega$, $A$ and $a$ are as in Theorem \ref{THM:well-posedness}. Let $u$ be the unique solution to \eqref{eqn:nonlinear_FMSE} for $g\in C^\infty_c(\Omega_e)$ with $\|g\|_{C^\infty_c(\Omega_e)}<\varepsilon$. Then for $\varepsilon>0$ small enough, we have
\[\|u\|_{H^s(\R^n)}\leq C\|g\|_{H^s(\R^n)}.\]  
\end{corollary}
\begin{proof}
In above proof of Theorem \ref{THM:well-posedness}, indeed we also have
\[\|u_0\|_{H^s(\R^n)}\leq C\|g\|_{H^s(\R^n)}\]
from \eqref{Cs estimate} and the well-posedness result of \cite{Covi2020} for the linear FMSE. Then using the $C^s$ bounds of $u_0$ and $v$ from the theorem, we obtain
\[\|\mathcal F(v)\|_{H^{-s}({\color{red}\Omega})}\leq \|\mathcal F(v)\|_{H^s({\color{red}\Omega})}\leq C(\varepsilon+\delta)(\|g\|_{H^s(\R^n)}+\|v\|_{H^s(\R^n)})\]
for $\delta, \varepsilon$ as in the proof of the theorem.
Then the well-posedness result of \cite{Covi2020} also implies the solution of \eqref{eqn:v} has $H^s$ bound, which satisfies
\[\|v\|_{H^s(\R^n)}\leq C\|\mathcal F(v)\|_{H^{-s}({\color{red}\Omega})}\leq C(\varepsilon+\delta)(\|g\|_{H^s(\R^n)}+\|v\|_{H^s(\R^n)}).\]
Lastly, by choosing the $\delta, \varepsilon$ pair in the proof of Theorem \ref{THM:well-posedness} to be further small to obtain
\[\|v\|_{H^s(\R^n)}\leq C\|g\|_{H^s(\R^n)}.\]
From $u=u_0+v$ and the estimates for $u_0$ and $v$ above, the proof is complete.
\end{proof}

\section{The inverse problem}\label{sec:inverse}
In this section, we will reconstruct the magnetic potential and nonlinear electric potential. 
We have showed in Theorem~\ref{THM:well-posedness} that for any 
$$g\in \mathcal{S}_{\varepsilon_0}:=\{g\in C_c^\infty(\Omega_e):\, \|g\|_{C^\infty_c(\Omega_e)}\leq \varepsilon_0\},$$ 
with $\varepsilon_0>0$ small enough,
there is a unique small solution $u_g\in C^s(\R^n)\cap H^s(\R^n)$ to the problem
\begin{align}\label{eqn:nonlinear_FMSE IP}
\begin{cases}
(-\Delta)^s_A u+a(x,u) =0 & \hbox{ in } \Omega,\\
u=g &  \hbox{ in } \Omega_e.\\
\end{cases} 
\end{align} 
 
To show that the map $g\mapsto u_g$ is differentiable in $\mathcal{S}_{\varepsilon_0}$, we consider 
for sufficiently small $\varepsilon>0$ and $f\in {C_c^\infty(\Omega_e)}$, let $u_{\varepsilon f}=u_{\varepsilon f}(x;\varepsilon)$ be the unique small solution to the problem
\begin{align}\label{eqn:nonlinear_FMSE 1}
\begin{cases}
(-\Delta)^s_A u+a(x,u) =0 & \hbox{ in } \Omega,\\
u=\varepsilon f  &  \hbox{ in } \Omega_e.\\
\end{cases} 
\end{align}

\subsection{Linearization}

We define the $k$-th derivative of the solution $u_{\varepsilon f}$ with respect to $\varepsilon$ by
$$u^{(k)}_{\varepsilon}(x) :=\frac{d^k}{d\varepsilon^k}u_{\varepsilon f}(x;\varepsilon)$$
for any positive integer $k$. We show that $u^{(k)}_\varepsilon$ satisfies various linear equations. 
{
\begin{lemma}\label{Lemma:differentibility} 
Let $\Omega$ be a bounded domain of $\R^n$ with smooth boundary. Assume $A\in C^s_c (\R^n\times\R^n;\R^n)$ and $a=a(x,z):\overline\Omega\times \R \to \R$ satisfy \eqref{eqn:MP_cond_A1}, \eqref{condition a} and \eqref{eqn:MP_cond_A2}.
Let $f\in {C_c^\infty(\Omega_e)}$, {$f\neq 0$} and $\varepsilon\in \R$ satisfy $|\varepsilon|<\frac{\varepsilon_0}{\|f\|_{C_c^\infty(\Omega_e)}}$ {for above $\varepsilon_0>0$}. 
Then we have $u_{\varepsilon f}$ is infinitely many times differentiable in $\varepsilon$ in $\left(-\frac{\varepsilon_0}{\|f\|_{C_c^\infty(\Omega_e)}}, \frac{\varepsilon_0}{\|f\|_{C_c^\infty(\Omega_e)}}\right)$. Moreover, we have for $k=1$, 
\[(-\Delta)_A^su^{(1)}_\varepsilon+\p_za(x,u_{\varepsilon f})u^{(1)}_\varepsilon=0\quad\textrm{ in }\Omega,\qquad u^{(1)}_\varepsilon=f\quad\textrm{ in }\Omega_e,\]
and for $k=2,3,\ldots$,
\begin{equation}\label{eqn:u_var_k}\left\{\begin{split}(-\Delta)^s_A u^{(k)}_\varepsilon+\p_za(x,u_{\varepsilon f})u^{(k)}_\varepsilon+\p_z^k a(x,u_{\varepsilon f}) (u_\varepsilon^{(1)})^k+R_{k-1}(a,u_{\varepsilon f})&=0\quad\textrm{ in }\Omega,\\
u^{(k)}_\varepsilon&=0\quad \textrm{ in }\Omega_e,
\end{split}\right.\end{equation}
where the term $R_{k-1}(a,u_{\varepsilon f})$ only involves ${\p_z^2a(x,u_{\varepsilon f}),\ldots,\p_z^{k-1} a(x,u_{\varepsilon f})}$ and $u^{(1)}_\varepsilon,\ldots,u^{(k-1)}_\varepsilon$. 
	
\end{lemma}
\begin{proof}
For $\varepsilon\in \left(-\frac{\varepsilon_0}{\|f\|_{C_c^\infty(\Omega_e)}}, \frac{\varepsilon_0}{\|f\|_{C_c^\infty(\Omega_e)}}\right)$ {and $\Delta\varepsilon\neq 0$}, 
set $$\widetilde{u}= {u_{(\varepsilon+\Delta\varepsilon)f} - u_{\varepsilon f}\over \Delta\varepsilon}.$$ By Taylor's formula, $\widetilde{u}$ is the solution to
	$$(-\Delta)^s_{A} \widetilde{u} + a^*(x)\widetilde{u}  = 0 \qquad \hbox{ in   $\Omega$}$$  with $\widetilde{u}=f$ in $\Omega_e$, where $$a^* (x) :=\int^1_0 \p_z a(x,su_{(\varepsilon+\Delta\varepsilon)f} + (1-s)u_{\varepsilon f}) \,ds$$ satisfies 
	$$(\nabla\cdot)^sA_{s||}+\int_{\R^n}|A(x,y)|^2~dy+a^*(x)\geq 0$$ for $x\in\Omega$ due to \eqref{eqn:MP_cond_A2}. By Lemma~\ref{lem:Cs_est},  
	\begin{equation*} \|\widetilde u\|_{C^s(\R^n)}\leq C\|f\|_{C_c^\infty(\Omega_e)},\end{equation*} 
	{which implies that
	\begin{align}\label{eqn:tildeu_Cs}
	\|u_{(\varepsilon+\Delta\varepsilon)f} - u_{\varepsilon f}\|_{C^s(\R^n)} \leq C|\Delta\varepsilon|\|f\|_{C_c^\infty(\Omega_e)}.
	\end{align}}
By {\eqref{eqn:tildeu_Cs}, we have}
	\begin{align*}
	\|a^*(x) -\p_z a(x,u_{\varepsilon f})\|_{L^\infty(\Omega)}  &= \|\int^1_0 (\p_z a(x,su_{(\varepsilon+\Delta\varepsilon)f} + (1-s)u_{\varepsilon f}) - \p_z a(x,u_{\varepsilon f})) \,ds \|_{L^\infty(\Omega)} \\
	&\leq C\|\p_z a(x,su_{(\varepsilon+\Delta\varepsilon)f} + (1-s)u_{\varepsilon f}) - \p_z a(x,u_{\varepsilon f})\|_{L^\infty(\Omega)}\\
	&\leq C\|\p_z^2 a\|_{L^\infty(\Omega)}\| u_{(\varepsilon+\Delta\varepsilon)f} - u_{\varepsilon f} \|_{C^s(\R^n)}\\
	&\leq C\|\p_z^2 a\|_{L^\infty(\Omega)} |\Delta\varepsilon|\|f\|_{C_c^\infty(\Omega_e)}.
	\end{align*}
This shows that as $\Delta\varepsilon\rightarrow0$, we have $\widetilde u\rightarrow u^{(1)}_\varepsilon$ in $C^s(\R^n)$, where $u^{(1)}_\varepsilon$ is the solution to
\[(-\Delta)_A^su^{(1)}_\varepsilon+\p_za(x, u_{\varepsilon f})u^{(1)}_\varepsilon=0\quad\textrm{ in }\Omega,\qquad u^{(1)}=f\quad\textrm{ in }\Omega_e.\]
In fact, if we denote 
	$
	w:=\widetilde{u}-u^{(1)}_\varepsilon,
	$
it is then a solution to 
	$$(-\Delta)^s_{A } w + a^*(x) w= (\p_z a(x,u_{\varepsilon f}) -a^*(x) )u^{(1)}_\varepsilon \quad\textrm{ in }\;\Omega,  \qquad w=0 \quad\textrm{ in }\; \Omega_e.$$
	Then we have the estimate
	\begin{align*} 
	\|w\|_{C^s(\R^n)} &\leq C\| (\p_z a(x,u_{\varepsilon f}) -a^*(x) )u^{(1)}_\varepsilon \|_{L^\infty(\Omega)} \\
	&\leq C\| \p_z a(x,u_{\varepsilon f}) -a^*(x)\|_{L^\infty(\Omega)} \|u^{(1)}_\varepsilon\|_{C^s(\R^n)} \\
	&\leq C\|\p_z^2 a\|_{L^\infty(\Omega)} |\Delta\varepsilon|\|f\|^2_{C_c^\infty(\Omega_e)}\rightarrow 0,\qquad\textrm{ as }\Delta\varepsilon\rightarrow0.
	\end{align*}
Similarly, set $\widetilde u^{(1)}:=\frac{u^{(1)}_{\varepsilon+\Delta\varepsilon}-u^{(1)}_\varepsilon}{\Delta\varepsilon}$. Then 
\[(-\Delta)_A^s\widetilde u^{(1)}+\frac1{\Delta\varepsilon}\left[\p_za\left(x,u_{(\varepsilon+\Delta\varepsilon)f}\right)u^{(1)}_{\varepsilon+\Delta\varepsilon}-\p_za\left(x,u_{\varepsilon f}\right)u^{(1)}_\varepsilon\right]=0,\]
giving
\[(-\Delta)_A^s\widetilde u^{(1)}+a^*_1(x)\widetilde uu^{(1)}_{\varepsilon+\Delta\varepsilon}+\p_za(x,u_{\varepsilon f})\widetilde u^{(1)}=0,\]
where
\[a^*_1(x):=\int_0^1\p_z^2a(x, su_{(\varepsilon+\Delta\varepsilon)f}+(1-s)u_{\varepsilon f})~ds\rightarrow \p_z^2a(x,u_{\varepsilon f})\]
in $L^\infty$ as $\Delta\varepsilon\rightarrow0$. 
Similar to above, as $\Delta\varepsilon\rightarrow 0$, we have $\widetilde u^{(1)}\rightarrow u^{(2)}_\varepsilon$ where
\[(-\Delta)_A^su^{(2)}_\varepsilon+\p_za(x,u_{\varepsilon f})u^{(2)}_\varepsilon+\p_z^2a(x, u_{\varepsilon f})(u^{(1)}_\varepsilon)^2=0.\]
{ Here we used the continuity of $u^{(1)}_\varepsilon$ in $\varepsilon$ which can be derived by following a similar argument in the derivation of \eqref{eqn:tildeu_Cs} above. } 
We apply the induction argument. Suppose \eqref{eqn:u_var_k} is true for index $k$ {and thus we have
\begin{align}\label{uell} 
	\widetilde u^{(\ell)}:= \frac{u^{(\ell)}_{\varepsilon+\Delta\varepsilon}-u^{(\ell)}_\varepsilon}{\Delta\varepsilon}\rightarrow u^{(\ell+1)}_\varepsilon
\end{align} 
	for $\ell=0,\ldots,k-1$. Here we denote $\widetilde u^{(0)} :=\widetilde u$.
} 
Set $\widetilde u^{(k)}:=\frac{u^{(k)}_{\varepsilon+\Delta\varepsilon}-u^{(k)}_\varepsilon}{\Delta\varepsilon}$. Then one can check for $k\geq 2$,
\begin{align}\label{eqn:uk}
\begin{split}(-\Delta)_A^s\widetilde u^{(k)}&+\frac{1}{\Delta\varepsilon}\left[\p_za(x,u_{(\varepsilon+\Delta\varepsilon)f})u_{\varepsilon+\Delta\varepsilon}^{(k)}-\p_za(x, u_{\varepsilon f})u_\varepsilon^{(k)}\right]\\
&+\frac1{\Delta\varepsilon}\left[\p_z^ka(x,u_{(\varepsilon+\Delta\varepsilon)f})(u^{(1)}_{\varepsilon+\Delta\varepsilon})^k-\p_z^ka(x, u_{\varepsilon f})(u^{(1)}_\varepsilon)^k\right]\\
&+\widetilde R_{k-1}=0, 
\end{split}
\end{align}
where the term $\widetilde R_{1}=0$ and $\widetilde R_{k-1}:=R_{k-1}(a, u_{(\varepsilon+\Delta\varepsilon)f})-R_{k-1}(a, u_{\varepsilon f})$ involves $a^*_2(x),\ldots, a^*_{k-1}(x)$, $\widetilde u, \widetilde u^{(1)},$ $\ldots, \widetilde u^{(k-1)}$, $\p_z^2a(x, u_{\varepsilon f}), \ldots, \p_z^{k-1}a(x,u_{\varepsilon f})$, $u_\varepsilon^{(1)},\ldots, u_\varepsilon^{(k-1)}$ and $u_{\varepsilon+\Delta\varepsilon}^{(1)},\ldots, u_{\varepsilon+\Delta\varepsilon}^{(k-1)}$ where 
\[a^*_l(x):=\int_0^1\p_z^{l+1}a(x, su_{(\varepsilon+\Delta\varepsilon)f}+(1-s)u_{\varepsilon f})~ds.\]
As $\Delta\varepsilon\rightarrow0$, {we can derive $\tilde u^{(k)} \rightarrow u^{(k+1)}_\varepsilon$ due to \eqref{condition a}, \eqref{uell} and a similar argument above. Here $u^{(k+1)}_\varepsilon$ satisfies}
\begin{align}\label{eqn:uk+1}
\begin{split}(-\Delta)_A^su^{(k+1)}_\varepsilon&+\p_za(x, u_{\varepsilon f})u^{(k+1)}_\varepsilon+\p_z^{k+1}a(x, u_{\varepsilon f})(u^{(1)}_\varepsilon)^{k+1}+R_k(a, u_{\varepsilon f})=0,
\end{split}
\end{align}
where 
\[R_k(a,u_{\varepsilon f})=\p_z^2a(x,u_{\varepsilon f})u^{(1)}_\varepsilon u^{(k)}_\varepsilon+k\p_z^ka(x,u_{\varepsilon f})(u^{(1)}_\varepsilon)^{k-1}u^{(2)}_\varepsilon+\lim_{\Delta\varepsilon\rightarrow0}\widetilde R_{k-1}\]
only involves ${\p_z^2a(x,u_{\varepsilon f}),\ldots,\p_z^{k} a(x,u_{\varepsilon f})}$ and $u^{(1)}_\varepsilon,\ldots,u^{(k)}_\varepsilon$. 
This completes the proof.
\end{proof}
}

\subsection{The DN map}

We now define the operator $B_{A,a}^s:H^s(\R^n) \times H^s(\R^n)\rightarrow \R $ by
$$
B_{A,a}^s [u,v] := \int_{\R^n}\int_{\R^n}\nabla^s_A u  \cdot\nabla^s_A v \,dydx + \int_{\Omega} a(x,u ) v dx.
$$
Now we give a definition of the DN map $\Lambda_{A,a}^s$.
\begin{lemma}\label{Lemma:DN-map} Let $\Omega$ be a bounded domain of $\R^n$ with smooth boundary. Assume $A\in C^s_c (\R^n\times\R^n;\R^n)$ and $a=a(x,z):\overline\Omega\times \R \to \R$ satisfy \eqref{eqn:MP_cond_A1}, \eqref{condition a} and \eqref{eqn:MP_cond_A2}. 
	There exists a bounded map $\Lambda_{A,a}^s : {X\cap\mathcal S_{\varepsilon_0}} \rightarrow X^*$ defined by
	$$
	   \langle \Lambda_{A,a}^s[g],[v]\rangle := B_{A,a}^s [u_{g} ,v ] \qquad \forall ~v\in H^s(\R^n),\;{g\in \mathcal S_{\varepsilon_0}},
	$$
	where $X$ is the quotient space $H^s(\R^n)\setminus \widetilde{H}^s(\Omega)$ and $u_g\in H^s(\R^n)$ solves $(-\Delta)^s_A u_g + a(x,u_g)=0$ in $\Omega$ with $u_g-g\in  \widetilde{H}^s(\Omega)$.
\end{lemma}
\begin{proof}
We first show that the definition of the DN map only depend on the equivalence classes. To this end, for any $\phi, \psi$ in $\widetilde{H}^s(\Omega)$, we have $u_{g+\phi}=u_g$ on $\R^n$ by uniqueness of the solution. Also, since $B_{A,s}^s$ is linear in the second component,  
\[B_{A,a}^s [u_{g+\phi} ,v+\psi ] =B_{A,a}^s[u_g,v+\psi]=B_{A,a}^s[u_g,v]+B_{A,a}^s[u_g,\psi],\]
where
\begin{align*}
 B_{A,a}^s [u_g,\psi]= \int_{\R^n}\int_{\R^n}\nabla^s_A u_{g} \cdot\nabla^s_A \psi \,dydx + \int_{\Omega} a(x,u_{g}) \psi dx=0
 \end{align*}
 by the weak formulation of the equation $(-\Delta)^s_Au_g+a(x,u_g)=0$ in $\Omega$ and $\psi\in \widetilde H^s(\Omega)$. 

Next the map is bounded because
\begin{align*}
|B_{A,a}^s [u_{g} ,v ]| 
&\leq C \|\nabla_A^s u_g\|_{L^2(\R^{2n})} \|\nabla_A^s v\|_{L^2(\R^{2n})} + \|a(x,u_g )\|_{L^2(\Omega)} \|v\|_{L^2(\R^n)}\\
&\leq C \|u_g\|_{H^s(\R^{n})} \|v\|_{H^s(\R^{n})}\\ 
& \leq C\|g\|_{H^s(\R^n)}\|v\|_{H^s(\R^n)},
\end{align*}
by Corollary \ref{Corollary:Hs_reg}.
Here we used the fact that $a(x,0)=0$ and Taylor's Theorem yield
\[\begin{split}
\|a(x,u_g)\|_{L^2({\Omega})}&= \left\|\left(\int_0^1\p_za(x,su_g(x))~ds\right)~u_g(x)\right\|_{L^2({\Omega})}\\
&\leq  \|\p_za(x,z)\|_{L^\infty(\Omega\times B_{C\varepsilon_0})}\|u_g\|_{L^2(\R^n)},
\end{split}\]
where $B_{C\varepsilon_0}$ stands for a ball with center at the origin and radius $C\varepsilon_0>0$.
This completes the proof.
\end{proof}

Lemma~\ref{Lemma:differentibility} implies that the solution $u_{\varepsilon f}=u_{\varepsilon f}(x;\varepsilon)$ to \eqref{eqn:nonlinear_FMSE 1} is differentiable with respect to $\varepsilon$ in the space $C^s(\R^n)$. To simplify the notation, we now denote the $k$-th derivative of the solution $u$ with respect to $\varepsilon$ at $\varepsilon=0$ by
\begin{equation}\label{eqn:def_u_k}
u^{(k)}(x)= {d^k\over d\varepsilon^k} \Big|_{\varepsilon=0} u_{\varepsilon f}(x;\varepsilon).
\end{equation}
Moreover, this allows us to take the $k$-th derivative $\Lambda^{(k),s}_{A,a}$ of the map $\Lambda_{A,a}^s$ with respect to $\varepsilon$ at $\varepsilon=0$:
$$
\langle \Lambda_{A,a}^{(k),s}[f],[v]\rangle := \lim_{\varepsilon\rightarrow 0} {d^k\over d\varepsilon^k} B_{A,a}^s [u_{\varepsilon f}(x;\varepsilon) ,v ], \qquad {f\in C_c^\infty(\Omega_e),\; v\in H^s(\R^n)}.
$$
\begin{lemma}\label{Lemma:DN_expansion}
Let $\Omega$ be a bounded domain of $\R^n$ with smooth boundary. Assume $A\in C^s_c(\R^n\times\R^n;\R^n)$ and $a=a(x,z):\overline\Omega\times \R \to \R$ satisfy \eqref{eqn:MP_cond_A1}, \eqref{condition a} and \eqref{eqn:MP_cond_A2}. Then we have 
\[\langle \Lambda_{A,a}^{(1),s}[f],[v]\rangle = \int_{\R^n}\int_{\R^n} \nabla^s_A  u^{(1)} \cdot\nabla^s_A v \,dydx + \int_{\Omega} \p_z a(x,0)u^{(1)}  v \,dx,\]
and for $k\geq 2$,
\[\begin{split}
	\langle \Lambda_{A,a}^{(k),s}[ f],[v]\rangle = & \int_{\R^n}\int_{\R^n} \nabla^s_A  u^{(k)} \cdot\nabla^s_A v \,dydx \\
	&+ \int_{\Omega} \LC {\p_za(x,0)u^{(k)}}+ \p_z^k a(x,0) (u^{(1)})^k+R_{k-1}(a,u) \RC  v \,dx,
\end{split}\]
for $v\in H^s(\R^n)$, where the term {$R_1=0$ and $R_{k-1}(a,u)$ only contains $\p_z^2 a(x,0),\ldots, \p_z^{k-1}a(x,0)$ and $u^{(1)},\ldots, u^{(k-1)}$}. 
\end{lemma}
\begin{proof}

	For any $v\in H^s(\R^n)$, {by the $C^s$ regularity of $A$ and $a$}, we can justify passing the limits into the integral to obtain 
	\begin{align*}
	&{d^k\over d\varepsilon^k} \Big|_{\varepsilon =0} B_{A,a}^s [u_{\varepsilon f}(x,\varepsilon),v ]\\
	&= \int_{\R^n}\int_{\R^n} \nabla^s_A  u^{(k)}(x)  \cdot\nabla^s_A v \,dydx + \int_{\Omega} \LC \p_z a(x,0)u^{(k)} + \ldots + \p_z^k a(x,0) (u^{(1)})^k\RC  v \,dx .
	\end{align*} 
\end{proof}

\subsection{Proof of Theorem~\ref{THH:uniqueness}}
Below are the steps to prove the Theorem \ref{THH:uniqueness}. Assume  
\begin{align}\label{eqn:DN}
	 \left.	\Lambda^s_{A_1,a_1}[g] \right|_{W_2} = 	  \left.	\Lambda^s_{A_2,a_2}[g] \right|_{W_2} \qquad \text{ for any }g\in {\mathcal S_{\varepsilon_0}\cap C^\infty_c(W_1)}.
\end{align} 
where $W_1,W_2$ are two arbitrary open subsets in $\Omega_e$.
For $f\in C_c^\infty(W_1)$ and $\varepsilon>0$ small enough, let $u_j$ be the solution to \eqref{eqn:LFMSE} with the exterior data $\varepsilon f$ and $A,a$ replaced by $A_j,a_j$ for $j=1,2$. The hypothesis \eqref{eqn:DN} yields
\[0=\langle (\Lambda_{A_1,a_1}^s-\Lambda^s_{A_2,a_2})[\varepsilon f],[v]\rangle=B_{A_1,a_1}^s[u_1(x,\varepsilon),v]-B_{A_2,a_2}^s[u_2(x,\varepsilon), v]\]
for any $v\in H^s(\R^n)$ satisfying $\supp (v\chi_{\Omega_e})\subset W_2$ where $\chi_{\Omega_e}$ denotes the characteristic function of $\Omega_e$. Immediately, we have the $k$-derivatives of $\Lambda^{s}_{A_j,a_j}$ satisfying
\[\left.\Lambda^{(k),s}_{A_1,a_1}[f] \right|_{W_2} = 	  \left.	\Lambda^{(k),s}_{A_2,a_2}[f] \right|_{W_2}\quad\textrm{ for }\; f\in C_c^\infty(W_1).\]

\subsubsection{The determination of the first order term}
We can recover $A$ and $\p_za$ up to a gauge. 
\begin{proposition}\label{Prop: recover A and a_1} 
Let $\Omega$ be a bounded domain of $\R^n$ with smooth boundary. Assume $A_j\in {\color{red}C^s_c}(\R^n\times\R^n;\R^n)$ with support in $\Omega\times\Omega$ and $a_j=a_j(x,z):\overline\Omega\times \R \to \R$ satisfy \eqref{eqn:MP_cond_A1}, \eqref{condition a} and \eqref{eqn:MP_cond_A2} for $j=1,2$. 
Suppose $\Lambda^s_{A_1, a_1}[g]|_{W_2}=\Lambda^s_{A_2,a_2}[g]|_{W_2}$ for $g\in \mathcal S_{\varepsilon_0}\cap C_c^\infty(W_1)$. Then we have 
\[(A_1, \p_za_1(\cdot,0))\sim(A_2, \p_za_2(\cdot,0)).\]
\end{proposition}
\begin{proof} 
	It is shown in \cite{Covi2020} that $\Lambda_{A,a}^{(1),s}$	
	is linear, symmetric and bounded on $X$, and $\Lambda_{A,a}^{(1),s}[f]|_{W_2}$ for all $f\in C^\infty_c(W_1)$ determines $\left(A(x,y), q(x)\right)$ up to the gauge $\sim$, where $q(x)=\p_za(x,0)$ here. 
\end{proof} 

We need the following lemma for determining the full nonlinear potential $a(x,u)$. 
\begin{lemma}\label{Lemma:uk}
Let $\Omega$ be a bounded domain of $\R^n$ with smooth boundary. Assume $A\in C^s_c(\R^n\times\R^n;\R^n)$ and $a=a(x,z):\overline\Omega\times \R \to \R$ satisfy \eqref{eqn:MP_cond_A1}, \eqref{condition a}, \eqref{eqn:MP_cond_A2} and 
\[(A_1,\p_za_1(\cdot,0))\sim(A_2,\p_za_2(\cdot,0)).\]
Given $f\in C_c^\infty(\Omega_e)$ and $\varepsilon>0$ small, let $u_j=u_j(x; \varepsilon)$ be the solution to 
\begin{align}\label{eqn:nonlinear_FMSE 2}
\begin{cases}
(-\Delta)^s_{A_j} u+a_j(x,u) =0 & \hbox{ in } \Omega,\\
u=\varepsilon f  &  \hbox{ in } \Omega_e,\\
\end{cases} 
\end{align}
for $j=1,2$.
Then we have $u_1^{(1)}=u_2^{(1)}$. 
Moreover, for $k\geq 2$, if 
\begin{align*} 
	\p_z^\ell a_1(x,0) = \p_z^\ell a_2(x,0)\;\textrm{ in }\Omega  \qquad \hbox{ for any $2\leq \ell\leq k$},
\end{align*}
then 
\begin{align}\label{The lower order terms}
    u^{(\ell)}_1 = u^{(\ell)}_2 \;\hbox{ in }\R^n  \qquad \hbox{ for any $2\leq \ell\leq k$},
\end{align} 
	where $u_j^{(\ell)}$ is defined as in \eqref{eqn:def_u_k}. 
\end{lemma}
\begin{proof}First, by Lemma \ref{Lemma:differentibility}, we have
{
$$
(-\Delta)_{A_j}^s u^{(1)}_j +\p_za_j(x,0)u^{(1)}_j =0\quad\textrm{ in }\Omega, \qquad u_j^{(1)}=f \quad\textrm{ in }\Omega_e. $$ 
Since $(A_1,\p_za_1(x,0))\sim(A_2,\p_za_2(x,0))$, then Lemma 3.3 and Lemma 3.8 in \cite{Covi2020}} and the well-posedness result yield that 
$$
u^{(1)}_1 = u^{(1)}_2\qquad\hbox{ in }\R^n.
$$

For $k=2$, by Lemma \ref{Lemma:differentibility} again, $u_j^{(2)}$ are solutions to
\[(-\Delta)_{A_j}^su_j^{(2)}+\p_za_j(x,0)u_j^{(2)}=-\p_z^2a_j(x,0)(u_j^{(1)})^2\]
When $\p_z^2a_1(x,0)=\p_z^2a_2(x,0)$ in $\Omega$, the equations for $u_j^{(2)}$, $j=1,2$ are identical and both $u_1^{(2)}$ and $u_2^{(2)}$ are zero in $\Omega_e$. Then we obtain $u_1^{(2)}=u_2^{(2)}$ by the well-posedness. 

Now suppose for $k\geq 2$, $\p_z^\ell a_1(x,0)=\p_z^\ell a_2(x,0)$ in $\Omega$ for $2\leq\ell\leq k+1$ and $u_1^{(\ell)}=u_2^{(\ell)}$ for $2\leq\ell\leq k$. The equations for $u_j^{(k+1)}$ is given by
\[(-\Delta)_{A_j}^su_j^{(k+1)}+\p_za_j(x,0)u_j^{(k+1)}=-\p_z^{k+1}a_j(x,0)(u_j^{(1)})^{k+1}+R_k(a_j, u_j).\]
Here $R_k(a_j, u_j)$ only involves $\p_z^2a_j(x,0),\ldots,\p_z^ka_j(x,0)$ and $u_j^{(1)},\ldots,u_j^{(k)}$ and thus we have $$R_k(a_1, u_1)=R_k(a_2, u_2).$$ By assumption, the equations for $u_1^{(k+1)}$ and $u_2^{(k+1)}$ are identical 
in $\Omega$ and both $u_1^{(k+1)}$ and $u_2^{(k+1)}$ are zero in $\Omega_e$, hence $u_1^{(k+1)}=u_2^{(k+1)}$ in $\R^n$. This completes the induction proof. 
\end{proof}

\noindent
\subsubsection{The proof of main result} Now we are ready to complete the proof of the main theorem by reconstructing the higher order term of Taylor expansion of $a$.
\begin{proof}[Proof of Theorem~\ref{THH:uniqueness}]

From Proposition \ref{Prop: recover A and a_1}, it remains to show $\p_z^k a_1(x,0) = \p_z^k a_2(x,0)$ for $x\in\Omega$ and $k\geq 2$ due to \eqref{condition a}.

First, let $v\in H^s(\R^n)$ be a function satisfying
\begin{equation}\label{eqn:linear_v}(-\Delta)_{A_1}^s v+\p_za_1(x,0)v=0\qquad\textrm{ in }\Omega\end{equation}
and $\supp(v\chi_{\Omega_e})\subset W_2$.
Note that here $(A_1, \p_za_1(\cdot,0))\sim(A_2,\p_za_2(\cdot,0))$ implies $(-\Delta)_{A_2}^s v+\p_za_2(x,0)v=0$ in $\Omega$. 
Then by Lemma \ref{Lemma:DN_expansion} and that $u_j^{(2)}=u_j^{(3)}=\cdots=0$ in $\Omega_e$, we have
\[\begin{split}
\langle\Lambda_{A_j,a_j}^{(k+1),s}[f],[v]\rangle&= \int_{\R^n}\int_{\R^n}\nabla_{A_j}^su_j^{(k+1)}\cdot\nabla_{A_j}^s v~dy~dx+\int_{\Omega}\p_za_j(x,0)u_j^{(k+1)}v~dx\\
&\quad+\int_{\Omega}\left[\p_z^{k+1}a_j(x,0)(u_j^{(1)})^{k+1}+R_k(a_j,u_j)\right]v~dx\\
&= \int_{\Omega}\left[\p_z^{k+1}a_j(x,0)(u_j^{(1)})^{k+1}+R_k(a_j,u_j)\right]v~dx,\qquad \textrm{ for } k\geq 1.
\end{split}
\]

By Lemma \ref{Lemma:uk}, we first have $u^{(1)}_1=u^{(1)}_2$. From $\Lambda_{A_1,a_1}^{(2),s}[f]|_{W_2}=\Lambda_{A_2,a_2}^{(2),s}[f]|_{W_2}$ for any $f\in C^\infty_c(W_1)$, we have 
\[0=\langle (\Lambda_{A_1,a_1}^{(2),s} - \Lambda_{A_2,a_2}^{(2),s})[ f],[v]\rangle=\int_{\Omega}(\p_z^2a_1(x,0)-\p_z^2a_2(x,0))(u^{(1)}_1)^2v~dx,\]
since $R_1(a_j,u_j)=0$.

By the Runge approximation property in Proposition \ref{prop:RAP}, 
for any $g\in L^2(\Omega)$, we can find a sequence of linear solutions $v_i\in H^s(\R^n)${, whose restriction to $\Omega_e$ is compactly supported in $W_2$, to \eqref{eqn:linear_v}} such that $v_i\rightarrow g$ in $L^2(\Omega)$. 
Plugging $v_i$ and let $i\rightarrow\infty$, we obtain 
\[\int_{\Omega}(\p_z^2a_1(x,0)-\p_z^2a_2(x,0))(u^{(1)}_1)^2g~dx=0\]
for any $g\in L^2(\Omega)$, which yields
\[(\p_z^2a_1(x,0)-\p_z^2a_2(x,0))(u^{(1)}_1)^2=0,\qquad x\in\Omega.\]
Now based on the strong maximum principle proved in Proposition~\ref{prop:MP_Strong}, we can choose suitable exterior data $u^{(1)}_1= f\geq 0$ and $f\not\equiv 0$ in $W_1$ so that $u_1^{(1)}>0$ in $\Omega$. Then
$(\p_z^{2}a_1(x,0)-\p_z^{2}a_2(x,0))$ must be zero in $\Omega$.

Now for any fixed positive integer $k>1$, we suppose that 
\begin{align}\label{lower order terms}
    \p_z^\ell a_1(x,0) = \p_z^\ell a_2(x,0) \quad\textrm{ in }\Omega,\qquad \textrm{ for }\; 2\leq \ell \leq k. 
\end{align} 
Then Lemma~\ref{Lemma:uk} yields that 
\begin{align}\label{lower order terms u}
    u^{(\ell)}_1 = u^{(\ell)}_2\quad \hbox{ in }\R^n,  \qquad \hbox{ for }\; 1\leq \ell \leq k. 
\end{align} 
Combining $\Lambda_{A_1,a_1}^{(k+1),s}[f]|_{W_2}=\Lambda_{A_2,a_2}^{(k+1),s}[f]|_{W_2}$ for $f\in C_c^\infty(W_1)$ with \eqref{lower order terms} and \eqref{lower order terms u}, one has
\[0=\int_{\R^n}(\p_z^{k+1}a_1(x,0)-\p_z^{k+1}a_2(x,0))(u_1^{(1)})^{k+1}v~dx\]
for $u^{(1)}_1$ and $v$ as above, which similarly proves $\p_z^{k+1}a_1(x,0)=\p_z^{k+1}a_2(x,0)$ in $\Omega$. By induction, this completes the proof.
\end{proof}

\noindent

\section{Acknowledgement}
R.-Y. Lai is partially supported by the National Science Foundation through grants DMS-2006731.

\bibliographystyle{abbrv}

\bibliography{FractionalMSE}

\begin{thebibliography}{10}

\bibitem{AZ2018}
Y.~Assylbekov and T.~Zhou.
\newblock Direct and inverse problems for the nonlinear time-harmonic maxwell
  equations in kerr-type media.
\newblock {\em To appear in Journal of Spectral Theory}, 2020.

\bibitem{AZ2020}
Y.~Assylbekov and T.~Zhou.
\newblock Inverse problems for nonlinear maxwell's equations with second
  harmonic generation.
\newblock {\em arXiv:2009.03467}, 2020.

\bibitem{BGU18}
S.~Bhattacharya, T.~Ghosh, and G.~Uhlmann.
\newblock Inverse problem for fractional-{L}aplacian with lower order non-local
  perturbations.
\newblock {\em arXiv preprint arXiv:1810.03567}, 2018.

\bibitem{BRSbook}
G.~B. Bisci, V.~D. Radulescu, and R.~Servadei.
\newblock Variational methods for nonlocal fractional problems.
\newblock {\em Encyclopedia of Mathematics and Its Applications Series}, 162,
  2016.

\bibitem{CNV2019}
C.~C\^arstea, G.~Nakamura, and M.~Vashisth.
\newblock Reconstruction for the coefficients of a quasilinear elliptic partial
  differential equation.
\newblock {\em Applied Mathematics Letters}, 98:121--127, 2019.

\bibitem{cekic2020calderon}
M.~Cekic, Y.-H. Lin, and A.~R{\"u}land.
\newblock The {C}alder{\'o}n problem for the fractional {S}chr{\"o}dinger
  equation with drift.
\newblock {\em Cal. Var. Partial Differential Equations},
  59(91):https://doi.org/10.1007/s00526--020--01740--6, 2020.

\bibitem{CLOP}
X.~Chen, M.~Lassas, L.~Oksanen, and G.~Paternain.
\newblock Detection of {H}ermitian connections in wave equations with cubic
  non-linearity.
\newblock {\em arXiv:1902.05711}, 2019.

\bibitem{Chung2014}
F.~Chung.
\newblock A partial data result for the magnetic {S}chr\"odinger inverse
  problem.
\newblock {\em Analysis and PDE}, 7:117--157, 2014.

\bibitem{ContTankovbook}
R.~Cont and P.~Tankov.
\newblock Financial modelling with jump processes.
\newblock {\em Chapman \& Hall/CRC Financial Mathematics Series}, 2004.

\bibitem{Covi2020}
G.~Covi.
\newblock An inverse problem for the fractional {S}chr\"odinger equation in a
  magnetic field.
\newblock {\em Inverse Problems}, 36:045004, 2020.

\bibitem{PietroSquassina2018}
P.~d'Avenia and M.~Squassina.
\newblock Ground states for fractional magnetic operators.
\newblock {\em ESAIM: Control, Optimisation and Calculus of Variations},
  24(1):1--24, 2018.

\bibitem{di2012hitchhiks}
E.~Di~Nezza, G.~Palatucci, and E.~Valdinoci.
\newblock Hitchhiker's guide to the fractional {S}obolev spaces.
\newblock {\em Bulletin des Sciences Math{\'e}matiques}, 136(5):521--573, 2012.

\bibitem{ER95}
G.~Eskin and J.~Ralston.
\newblock Inverse scattering problem for the {S}chr\"odinger equation with
  magnetic potential at fixed energy.
\newblock {\em Comm. Math. Phys.}, 173:199--224, 1995.

\bibitem{FO2019}
A.~Feizmohammadi and L.~Oksanen.
\newblock An inverse problem for a semi-linear elliptic equation in
  {R}iemannian geometries.
\newblock {\em Journal of Differential Equations}, 269:4683--4719, 2020.

\bibitem{Dos2009}
D.~D.~S. Ferreira, C.~Kenig, J.~Sj\"ostrand, and G.~Uhlmann.
\newblock Determining a magnetic {S}chr\"odinger operator from partial {C}auchy
  data.
\newblock {\em Comm. Math. Phys.}, 271(2):467--488, 2009.

\bibitem{ghosh2017calder}
T.~Ghosh, Y.-H. Lin, and J.~Xiao.
\newblock The {C}alder\'{o}n problem for variable coefficients nonlocal
  elliptic operators.
\newblock {\em Communications in Partial Differential Equations},
  42(12):1923--1961, 2017.

\bibitem{GRSU18}
T.~Ghosh, A.~R{\"u}land, M.~Salo, and G.~Uhlmann.
\newblock Uniqueness and reconstruction for the fractional {C}alder{\'o}n
  problem with a single measurement.
\newblock {\em Journal of Functional Analysis}, page 108505, 2020.

\bibitem{ghosh2016calder}
T.~Ghosh, M.~Salo, and G.~Uhlmann.
\newblock The {C}alder{\'o}n problem for the fractional {S}chr{\"o}dinger
  equation.
\newblock {\em Analysis \& PDE}, 13(2):455--475, 2020.

\bibitem{GT2011}
C.~Guillarmou and L.~Tzou.
\newblock Identification of a connection from {C}auchy data on a {R}iemann
  surface with bounday.
\newblock {\em Geom. Funct. Ana}, 21:393--418, 2011.

\bibitem{Haberman16}
B.~Haberman.
\newblock Unique determination of a magnetic {S}chr\"odinger operator with
  unbounded magnetic potential from boundary data.
\newblock {\em Int. Math. Res. Not.}, 4:1080--1128, 2018.

\bibitem{harrach2017nonlocal-monotonicity}
B.~Harrach and Y.-H. Lin.
\newblock Monotonicity-based inversion of the fractional {S}chr\"odinger
  equation {I}. {P}ositive potentials.
\newblock {\em SIAM Journal on Mathematical Analysis}, 51(4):3092--3111, 2019.

\bibitem{harrach2020monotonicity}
B.~Harrach and Y.-H. Lin.
\newblock Monotonicity-based inversion of the fractional {S}chr\"odinger
  equation {II}. {G}eneral potentials and stability.
\newblock {\em SIAM Journal on Mathematical Analysis}, 52(1):402--436, 2020.

\bibitem{Sun2002}
D.~Hervas and Z.~Sun.
\newblock An inverse boundary value problem for quasilinear elliptic equations.
\newblock {\em Communications in Partial Differential Equations},
  27:2449--2490, 2002.

\bibitem{IUY2012}
O.~Imanuvilov, G.~Uhlmann, and M.~Yamamoto.
\newblock Partial {C}auchy data for general second order elliptic operators in
  two dimensions.
\newblock {\em Publ. Res. Inst. Math. Sci.}, 48:971--1055, 2012.

\bibitem{Isakov93}
V.~Isakov.
\newblock On uniqueness in inverse problems for semilinear parabolic equations.
\newblock {\em Archive for Rational Mechanics and Analysis}, 124(1):1--12,
  1993.

\bibitem{victor01}
V.~Isakov.
\newblock Uniqueness of recovery of some quasilinear partial differential
  equations.
\newblock {\em Commun. in partial differential equations}, 26(11,
  12):1947--1973, 2001.

\bibitem{victorN}
V.~Isakov and A.~Nachman.
\newblock Global uniqueness for a two-dimensional elliptic inverse problem.
\newblock {\em Trans.of AMS}, 347:3375--3391, 1995.

\bibitem{isakov1994global}
V.~Isakov and J.~Sylvester.
\newblock Global uniqueness for a semilinear elliptic inverse problem.
\newblock {\em Communications on Pure and Applied Mathematics},
  47(10):1403--1410, 1994.

\bibitem{Kang2002}
H.~Kang and G.~Nakamura.
\newblock Identification of nonlinearity in a conductivity equation via the
  {D}irichlet-to-{N}eumann map.
\newblock {\em Inverse Problems}, 18:1079--1088, 2002.

\bibitem{KU14}
K.~Krupchyk and G.~Uhlmann.
\newblock Uniqueness in an inverse boundary problem for a magnetic
  {S}chr\"odinger operator with a bounded magnetic potential.
\newblock {\em Comm. Math. Phys.}, 327:993--1009, 2014.

\bibitem{KU201909}
K.~Krupchyk and G.~Uhlmann.
\newblock Partial data inverse problems for semilinear elliptic equations with
  gradient nonlinearities.
\newblock {\em To appear in Math. Res. Lett.}, 2019.

\bibitem{KU201905}
K.~Krupchyk and G.~Uhlmann.
\newblock A remark on partial data inverse problems for semilinear elliptic
  equations.
\newblock {\em Proceedings of the AMS}, 148:681--685, 2020.

\bibitem{KLU2018}
Y.~Kurylev, M.~Lassas, and G.~Uhlmann.
\newblock Inverse problems for {L}orentzian manifolds and non-linear hyperbolic
  equations.
\newblock {\em Invent. Math.}, 212(3):781--857, 2018.

\bibitem{Lai2011}
R.-Y. Lai.
\newblock Global uniqueness for an inverse problem for the magnetic
  {S}chr\"odinger operator.
\newblock {\em Inverse Problems and Imaging}, 5:59--74, 2011.

\bibitem{lai2019global}
R.-Y. Lai and Y.-H. Lin.
\newblock Global uniqueness for the fractional semilinear {S}chr{\"o}dinger
  equation.
\newblock {\em Proc. Amer. Math. Soc.}, 147(3):1189--1199, 2019.

\bibitem{LaiL2020}
R.-Y. Lai and Y.-H. Lin.
\newblock Inverse problems for fractional semilinear elliptic equations.
\newblock {\em arXiv:2004.00549}, 2020.

\bibitem{LLR2019calder}
R.-Y. Lai, Y.-H. Lin, and A.~R{\"u}land.
\newblock The {C}alder\'on problem for a space-time fractional parabolic
  equation.
\newblock {\em SIAM Journal on Mathematical Analysis}, 52:2655--2688, 2020.

\bibitem{LaiOhm20}
R.-Y. Lai and L.~Ohm.
\newblock Inverse problems for the fractional {L}aplace equation with lower
  order nonlinear perturbations.
\newblock {\em arXiv:2009.07883}, 2020.

\bibitem{LUY2020}
R.-Y. Lai, G.~Uhlmann, and Y.~Yang.
\newblock Reconstruction of the collision kernel in the nonlinear {B}oltzmann
  equation.
\newblock {\em To appear in SIAM J. Math. Anal.}, 2020.

\bibitem{LaiZhou2020}
R.-Y. Lai and T.~Zhou.
\newblock Partial data inverse problems for nonlinear magnetic {S}chr\"odinger
  equations.
\newblock {\em arXiv:2007.02475}, 2020.

\bibitem{LLLS201905}
M.~Lassas, T.~Liimatainen, Y.-H. Lin, and M.~Salo.
\newblock Partial data inverse problems and simultaneous recovery of boundary
  and coefficients for semilinear elliptic equations.
\newblock {\em To appear in Rev. Mat. Iberoam.}, 2019.

\bibitem{LLLS201903}
M.~Lassas, T.~Liimatainen, Y.-H. Lin, and M.~Salo.
\newblock Inverse problems for elliptic equations with power type
  nonlinearities.
\newblock {\em J. Math. Pures Appl.}, 145:44--82, 2021.

\bibitem{LLPT2020_1}
M.~Lassas, T.~Liimatainen, L.~Potenciano-Machado, and T.~Tyni.
\newblock Uniqueness and stability of an inverse problem for a semi-linear wave
  equation.
\newblock {\em arXiv:2006.13193}, 2020.

\bibitem{LUW2018}
M.~Lassas, G.~Uhlmann, and Y.~Wang.
\newblock Inverse problems for semilinear wave equations on {L}orentzian
  manifolds.
\newblock {\em Comm. Math. Phys.}, 360(2):555--609, 2018.

\bibitem{Li2019}
L.~Li.
\newblock The {C}alder\'on problem for the fractional magnetic operator.
\newblock {\em Inverse Problems}, 36:075003, 2020.

\bibitem{Li2020}
L.~Li.
\newblock A semilinear inverse problem for the fractional magnetic laplacian.
\newblock {\em arXiv:2005.06714}, 2020.

\bibitem{Lin202004}
Y.-H. Lin.
\newblock Monotonicity-based inversion of fractional semilinear elliptic
  equations with power type nonlinearities.
\newblock {\em arXiv:2005.07163}, 2020.

\bibitem{NSU95}
G.~Nakamura, Z.~Sun, and G.~Uhlmann.
\newblock Global identifiability for an inverse problem for the {S}chr\"odinger
  equation in a magnetic field.
\newblock {\em Math. Ann.}, 303:377--388, 1995.

\bibitem{PSU10}
L.~P\"aiv\"arinta, M.~Salo, and G.~Uhlmann.
\newblock Inverse scattering for the magnetic {S}chr\"odinger operator.
\newblock {\em J. Funct. Anal.}, 259:1771--1798, 2010.

\bibitem{ros2014dirichlet}
X.~Ros-Oton and J.~Serra.
\newblock The {D}irichlet problem for the fractional {L}aplacian: regularity up
  to the boundary.
\newblock {\em Journal de Math{\'e}matiques Pures et Appliqu{\'e}es},
  101(3):275--302, 2014.

\bibitem{ruland2017fractional}
A.~R{\"u}land and M.~Salo.
\newblock The fractional {C}alder\'{o}n problem: low regularity and stability.
\newblock {\em Nonlinear Analysis}, 193:111529, 2020.

\bibitem{Sun96}
Z.~Sun.
\newblock On a quasilinear inverse boundary value problem.
\newblock {\em Math. Z.}, 221(2):293--305, 1996.

\bibitem{SunUhlmann97}
Z.~Sun and G.~Uhlmann.
\newblock Inverse problems in quasilinear anisotropic media.
\newblock {\em American Journal of Mathematics}, 119(4):771--797, 1997.

\end{thebibliography}

\end{document}